\documentclass[10pt]{amsart}
\usepackage[dvips]{graphicx}

\usepackage{mathrsfs}
\usepackage{amssymb}
\usepackage{amsfonts}
\usepackage{latexsym}
\usepackage{amsthm}

\usepackage{graphicx}
\def\lb{\label}

\newcommand{\er}[1]{\textrm{(\ref{#1})}}

\begin{document}


\renewcommand{\theequation}{\arabic{section}.\arabic{equation}}
\theoremstyle{plain}
\newtheorem{theorem}{\bf Theorem}[section]
\newtheorem{lemma}[theorem]{\bf Lemma}
\newtheorem{corollary}[theorem]{\bf Corollary}
\newtheorem{proposition}[theorem]{\bf Proposition}
\newtheorem{definition}[theorem]{\bf Definition}
\newtheorem{remark}[theorem]{\it Remark}

\def\a{\alpha}  \def\cA{{\mathcal A}}     \def\bA{{\bf A}}  \def\mA{{\mathscr A}}
\def\b{\beta}   \def\cB{{\mathcal B}}     \def\bB{{\bf B}}  \def\mB{{\mathscr B}}
\def\g{\gamma}  \def\cC{{\mathcal C}}     \def\bC{{\bf C}}  \def\mC{{\mathscr C}}
\def\G{\Gamma}  \def\cD{{\mathcal D}}     \def\bD{{\bf D}}  \def\mD{{\mathscr D}}
\def\d{\delta}  \def\cE{{\mathcal E}}     \def\bE{{\bf E}}  \def\mE{{\mathscr E}}
\def\D{\Delta}  \def\cF{{\mathcal F}}     \def\bF{{\bf F}}  \def\mF{{\mathscr F}}
\def\c{\chi}    \def\cG{{\mathcal G}}     \def\bG{{\bf G}}  \def\mG{{\mathscr G}}
\def\z{\zeta}   \def\cH{{\mathcal H}}     \def\bH{{\bf H}}  \def\mH{{\mathscr H}}
\def\e{\eta}    \def\cI{{\mathcal I}}     \def\bI{{\bf I}}  \def\mI{{\mathscr I}}
\def\p{\psi}    \def\cJ{{\mathcal J}}     \def\bJ{{\bf J}}  \def\mJ{{\mathscr J}}
\def\vT{\Theta} \def\cK{{\mathcal K}}     \def\bK{{\bf K}}  \def\mK{{\mathscr K}}
\def\k{\kappa}  \def\cL{{\mathcal L}}     \def\bL{{\bf L}}  \def\mL{{\mathscr L}}
\def\l{\lambda} \def\cM{{\mathcal M}}     \def\bM{{\bf M}}  \def\mM{{\mathscr M}}
\def\L{\Lambda} \def\cN{{\mathcal N}}     \def\bN{{\bf N}}  \def\mN{{\mathscr N}}
\def\m{\mu}     \def\cO{{\mathcal O}}     \def\bO{{\bf O}}  \def\mO{{\mathscr O}}
\def\n{\nu}     \def\cP{{\mathcal P}}     \def\bP{{\bf P}}  \def\mP{{\mathscr P}}
\def\r{\varrho} \def\cQ{{\mathcal Q}}     \def\bQ{{\bf Q}}  \def\mQ{{\mathscr Q}}
\def\s{\sigma}  \def\cR{{\mathcal R}}     \def\bR{{\bf R}}  \def\mR{{\mathscr R}}
\def\S{\Sigma}  \def\cS{{\mathcal S}}     \def\bS{{\bf S}}  \def\mS{{\mathscr S}}
\def\t{\tau}    \def\cT{{\mathcal T}}     \def\bT{{\bf T}}  \def\mT{{\mathscr T}}
\def\f{\phi}    \def\cU{{\mathcal U}}     \def\bU{{\bf U}}  \def\mU{{\mathscr U}}
\def\F{\Phi}    \def\cV{{\mathcal V}}     \def\bV{{\bf V}}  \def\mV{{\mathscr V}}
\def\P{\Psi}    \def\cW{{\mathcal W}}     \def\bW{{\bf W}}  \def\mW{{\mathscr W}}
\def\o{\omega}  \def\cX{{\mathcal X}}     \def\bX{{\bf X}}  \def\mX{{\mathscr X}}
\def\x{\xi}     \def\cY{{\mathcal Y}}     \def\bY{{\bf Y}}  \def\mY{{\mathscr Y}}
\def\X{\Xi}     \def\cZ{{\mathcal Z}}     \def\bZ{{\bf Z}}  \def\mZ{{\mathscr Z}}
\def\O{\Omega}

\newcommand{\mc}{\mathscr {c}}

\newcommand{\gA}{\mathfrak{A}}          \newcommand{\ga}{\mathfrak{a}}
\newcommand{\gB}{\mathfrak{B}}          \newcommand{\gb}{\mathfrak{b}}
\newcommand{\gC}{\mathfrak{C}}          \newcommand{\gc}{\mathfrak{c}}
\newcommand{\gD}{\mathfrak{D}}          \newcommand{\gd}{\mathfrak{d}}
\newcommand{\gE}{\mathfrak{E}}
\newcommand{\gF}{\mathfrak{F}}           \newcommand{\gf}{\mathfrak{f}}
\newcommand{\gG}{\mathfrak{G}}           
\newcommand{\gH}{\mathfrak{H}}           \newcommand{\gh}{\mathfrak{h}}
\newcommand{\gI}{\mathfrak{I}}           \newcommand{\gi}{\mathfrak{i}}
\newcommand{\gJ}{\mathfrak{J}}           \newcommand{\gj}{\mathfrak{j}}
\newcommand{\gK}{\mathfrak{K}}            \newcommand{\gk}{\mathfrak{k}}
\newcommand{\gL}{\mathfrak{L}}            \newcommand{\gl}{\mathfrak{l}}
\newcommand{\gM}{\mathfrak{M}}            \newcommand{\gm}{\mathfrak{m}}
\newcommand{\gN}{\mathfrak{N}}            \newcommand{\gn}{\mathfrak{n}}
\newcommand{\gO}{\mathfrak{O}}
\newcommand{\gP}{\mathfrak{P}}             \newcommand{\gp}{\mathfrak{p}}
\newcommand{\gQ}{\mathfrak{Q}}             \newcommand{\gq}{\mathfrak{q}}
\newcommand{\gR}{\mathfrak{R}}             \newcommand{\gr}{\mathfrak{r}}
\newcommand{\gS}{\mathfrak{S}}              \newcommand{\gs}{\mathfrak{s}}
\newcommand{\gT}{\mathfrak{T}}             \newcommand{\gt}{\mathfrak{t}}
\newcommand{\gU}{\mathfrak{U}}             \newcommand{\gu}{\mathfrak{u}}
\newcommand{\gV}{\mathfrak{V}}             \newcommand{\gv}{\mathfrak{v}}
\newcommand{\gW}{\mathfrak{W}}             \newcommand{\gw}{\mathfrak{w}}
\newcommand{\gX}{\mathfrak{X}}               \newcommand{\gx}{\mathfrak{x}}
\newcommand{\gY}{\mathfrak{Y}}              \newcommand{\gy}{\mathfrak{y}}
\newcommand{\gZ}{\mathfrak{Z}}             \newcommand{\gz}{\mathfrak{z}}

\def\ve{\varepsilon}   \def\vt{\vartheta}    \def\vp{\varphi}    \def\vk{\varkappa}

\def\A{{\mathbb A}} \def\B{{\mathbb B}} \def\C{{\mathbb C}}
\def\dD{{\mathbb D}} \def\E{{\mathbb E}} \def\dF{{\mathbb F}} \def\dG{{\mathbb G}} \def\H{{\mathbb H}}\def\I{{\mathbb I}} \def\J{{\mathbb J}} \def\K{{\mathbb K}} \def\dL{{\mathbb L}}\def\M{{\mathbb M}} \def\N{{\mathbb N}} \def\O{{\mathbb O}} \def\dP{{\mathbb P}} \def\R{{\mathbb R}} \def\dQ{{\mathbb Q}}
\def\S{{\mathbb S}} \def\T{{\mathbb T}} \def\U{{\mathbb U}} \def\V{{\mathbb V}}\def\W{{\mathbb W}} \def\X{{\mathbb X}} \def\Y{{\mathbb Y}} \def\Z{{\mathbb Z}}

\newcommand{\1}{\mathbbm 1}
\newcommand{\dd}    {\, \mathrm d}



\def\la{\leftarrow}              \def\ra{\rightarrow}            \def\Ra{\Rightarrow}
\def\ua{\uparrow}                \def\da{\downarrow}
\def\lra{\leftrightarrow}        \def\Lra{\Leftrightarrow}


\def\lt{\biggl}                  \def\rt{\biggr}
\def\ol{\overline}               \def\wt{\widetilde}
\def\no{\noindent}


\let\ge\geqslant                 \let\le\leqslant
\def\lan{\langle}                \def\ran{\rangle}
\def\/{\over}                    \def\iy{\infty}
\def\sm{\setminus}               \def\es{\emptyset}
\def\ss{\subset}                 \def\ts{\times}
\def\pa{\partial}                \def\os{\oplus}
\def\om{\ominus}                 \def\ev{\equiv}
\def\iint{\int\!\!\!\int}        \def\iintt{\mathop{\int\!\!\int\!\!\dots\!\!\int}\limits}
\def\el2{\ell^{\,2}}             \def\1{1\!\!1}
\def\sh{\sharp}
\def\wh{\widehat}

\def\all{\mathop{\mathrm{all}}\nolimits}
\def\where{\mathop{\mathrm{where}}\nolimits}
\def\as{\mathop{\mathrm{as}}\nolimits}
\def\Area{\mathop{\mathrm{Area}}\nolimits}
\def\arg{\mathop{\mathrm{arg}}\nolimits}
\def\const{\mathop{\mathrm{const}}\nolimits}
\def\det{\mathop{\mathrm{det}}\nolimits}
\def\diag{\mathop{\mathrm{diag}}\nolimits}
\def\diam{\mathop{\mathrm{diam}}\nolimits}
\def\dim{\mathop{\mathrm{dim}}\nolimits}
\def\dist{\mathop{\mathrm{dist}}\nolimits}
\def\Im{\mathop{\mathrm{Im}}\nolimits}
\def\Iso{\mathop{\mathrm{Iso}}\nolimits}
\def\Ker{\mathop{\mathrm{Ker}}\nolimits}
\def\Lip{\mathop{\mathrm{Lip}}\nolimits}
\def\rank{\mathop{\mathrm{rank}}\limits}
\def\Ran{\mathop{\mathrm{Ran}}\nolimits}
\def\Re{\mathop{\mathrm{Re}}\nolimits}
\def\Res{\mathop{\mathrm{Res}}\nolimits}
\def\res{\mathop{\mathrm{res}}\limits}
\def\sign{\mathop{\mathrm{sign}}\nolimits}
\def\span{\mathop{\mathrm{span}}\nolimits}
\def\supp{\mathop{\mathrm{supp}}\nolimits}
\def\Tr{\mathop{\mathrm{Tr}}\nolimits}
\def\BBox{\hspace{1mm}\vrule height6pt width5.5pt depth0pt \hspace{6pt}}


\newcommand\nh[2]{\widehat{#1}\vphantom{#1}^{(#2)}}
\def\dia{\diamond}

\def\Oplus{\bigoplus\nolimits}




\def\qqq{\qquad}
\def\qq{\quad}
\let\ge\geqslant
\let\le\leqslant
\let\geq\geqslant
\let\leq\leqslant

\newcommand{\ca}{\begin{cases}}
\newcommand{\ac}{\end{cases}}
\newcommand{\ma}{\begin{pmatrix}}
\newcommand{\am}{\end{pmatrix}}
\renewcommand{\[}{\begin{equation}}
\renewcommand{\]}{\end{equation}}
\def\bu{\bullet}

\title[{Global transformations
preserving  spectral data }]
{Global transformations
preserving  Sturm-Liouville spectral data }

\date{\today}

\author[Hiroshi Isozaki]{Hiroshi Isozaki}
\address{Institute of Mathematics,
University of Tsukuba,
Tsukuba, 305-8571, Japan,\\
\ isozakih@math.tsukuba.ac.jp}
\author[Evgeny L. Korotyaev]{Evgeny L. Korotyaev}
\address{
Mathematical Physics Department, Faculty of Physics,
St. Petersburg State University,  Ulianovskaya 2, St. Petersburg, 198904, Russia,
 \ korotyaev@gmail.com}


\begin{abstract}
We show the existence of a real  analytic isomorphism between a space of impedance function $\rho$ of the Sturm-Liouville problem $- \rho^{-2}\big(\rho^2f')' + uf$ on $(0,1)$,
where $u$ is a function of $\rho, \rho', \rho''$, and
 that of potential $p$ of the Schr{\"o}dinger equation $- y'' + py$ on $(0,1)$, keeping their boundary conditions and spectral data. This mapping is  associated with the classical Liouville transformation $f \to \rho f$, and yields a global isomorphism between solutions to  inverse problems for the Sturm-Liouville equations of the impedance form and those to the Schr{\"o}dinger equations.
\end{abstract}

\maketitle

\section {Introduction and main results}
\setcounter{equation}{0}

\subsection{Louville transformation}
It is well-known that the Liouville transformation $f \to \rho f$
maps the solution to the Sturm-Liouville equation
\begin{equation}
- \rho^{-2}\big(\rho^2 f'\big)' + uf = \lambda f,
\label{IntroEquation}
\end{equation}
called the equation of {\it impedance form}, to that for  the
Schr{\"o}dinger equation
\begin{equation}
- y'' + py = \lambda y,
\label{IntroSchrod}
\end{equation}
and reduces the issues for the former to those for the latter. A
characteristic feature of the 1-dimensional inverse problem for the
Schr{\"o}dinger equation is that we know the global structure of
solutions, i.e. the existence of a real analytic isomorphism between
certain Hilbert space of potentials and that for the spectral data.
Therefore, one may well think of transforming back this isomorphism
by the inverse Liouville transformation to solve the inverse problem
for the Sturm-Liouville equation of impedance form. However, it is
by no means obvious to find a proper form of the associated
transformation between two inverse problems, in particular,
the function spaces for the impedance functions  and the potentials.
Thus  there has been no general result so far  except for some
restricted cases. In this paper, we show the existence of an
analytic isomorphism between a space of impedance functions $\rho$
and that for the potentials $p$, which also preserves the boundary
conditions and the spectral data, i.e. eigenvalues and norming
constants. Therefore, it gives an isomorphism between the solutions
of inverse problems for these equations. Since the inverse problem
for (\ref{IntroSchrod}) has already been solved, we can thus solve
the inverse problem for (\ref{IntroEquation}), in particular, the
characterization problem. This result can be applied to the
Laplacian on the surface of revolution \cite{G95} and to the
propagation of wave in the media periodic in radius \cite{TS09},
\cite{TS10}.  We shall discuss these problems elsewhere \cite{IK13}.

Suppose we are given a Sturm-Liouville operator  $-\D_{q,u}$ defined in
  $L^2((0,1);\r^2dx)$, where $\r=\r(x)>0$, having the form
\[
\lb{if1}
-\D_{q,u} f=-{1\/ \r^2}(\r^2f')'+u f,\qqq\qqq q={\r\ '\/\r},
\]
equipped with the boundary condition
\[
\lb{if2}
f'(0)-af(0)=0,\qquad f'(1)+b f(1)=0,\qquad a,b\in \R\cup \{\iy\}.
\]
In view of (\ref{if1}), we take $\r(x)$ as follows
\[
\lb{dro}
\r(x)=e^{Q(x)},\qqq Q=\int_0^xq(t)dt.
\]
Using the unitary transformation $\mU$ defined by
\[
\begin{aligned}
\lb{dU}
 \mU: L^2((0,1);\r^2 dx)\to L^2((0,1);dx),\\
 f\to \mU f= \r f,\qqq \r=e^{\int_0^xq(t)dt},
\end{aligned}
\]
we can transform the operator $-\D_{q,u}$ into the Schr\"odinger operator on
$L^2((0,1);dx)$:
\[
\lb{5x} \mU (-\D_q) \mU^{-1}=
 -\r^{-1}\frac{d}{dx} \Big(\r^2\frac{d}{dx} \r^{-1}\Big)+u=
 - \frac{d^2}{dx^2} + q' + q^2 + u.
\]
We now put
\[
\begin{aligned}
\lb{6x} S_p =-{d^2\/dx^2}+p,\qqq
\end{aligned}
\]
\begin{equation}
\lb{dP1}
\begin{aligned}
\ca  & p=q'+q^2+u-c_0, \qqq  \\
&{\displaystyle  c_0=\int_0^1 (q^2+u)dt},\ac
\end{aligned}
\end{equation}
and denote $L^2((0,1);dx)$  by $L^2(0,1)$.


\subsection{First main result - Analytic isomorphism}
Let $H^m$ be the Sobolev space of order $m$ on $(0,1)$ :
$$
H^m = \left\{ q \in L^2(0,1) \, ; \, q^{(k)} \in L^2(0,1),
 \ 0 \leq k \leq m\right\},
$$
and introduce the following space of real functions
\begin{equation}
\begin{split}
&\mW_1^0 =\rt\{q \in H^1\, ;\, q(0)=q(1)=0\rt\}, \qqq
\\
&\mH_0 =\rt\{q\in L^2(0,1)\, ; \  \int
_0^1q(x)dx=0\rt\}, \\
& \mH_\a = \mH_0\cap H^{\alpha},\qqq \a\ge 0,\\
\end{split}
\end{equation}
equipped with the norms
$$
\|q\|^2_{\mW_1^0}=\|q'\|^2=\int_0^1|q'(x)|^2dx,\qqq
\|q\|^2_{\mH_\a}=\|q^{(\a)}\|^2=\int_0^1|q^{(\a)}(x)|^2dx.
$$
Here and what follows, $\|\cdot\|$ denotes the norm of $L^2(0,1)$.
To show that they define the norms, we have only to pass to the
Fourier series. Define the space $L_{even}^2(0,1)$ of even functions
and the space  $L_{odd}^2(0,1)$ of odd functions by
\[
\begin{split}
\lb{oeL}
L_{odd}^2(0,1) &=\rt\{q\in L^2(0,1): q(x)=-q(1-x), \qq \forall
 \ x\in (0,1)\rt\},\\
L_{even}^2(0,1) &=\rt\{q\in L^2(0,1): q(x)=q(1-x), \qq \forall
 \ x\in (0,1)\rt\},
\end{split}
\]
and for $\o=even$ or  $\o=odd$ we put
\[
\begin{aligned}
\lb{oe}
\mW_1^{0,\o}= \mW_1^0\cap L_{\o}^2(0,1),\qqq
\mH_\a^{\o} = \mH_\a\cap L_{\o}^2(0,1), \qq \a\ge 0.
\end{aligned}
\]

We also introduce the space $\ell^2_{\a}$ of real sequences $h=(h_n)_1^{\iy }$,
equipped with the norm
\[
\lb{1.0}
\|h\|_{\a}^2=2\sum _{n\ge 1}(2\pi n)^{2\a}|h_n|^2,\qqq  \a\in \R,
\]
and  let $\ell^2=\ell_0^2$.

We assume that the potential $u$ is related to $q$ in the following way.

\bigskip
\noindent
\no  {\bf Condition U.} {\it
The potential $u$ has the form
\[
u=u_1(q)+u_2(Q),
\]
$Q$ being defined in (\ref{dro}),
with the following properties. \\
\noindent
(1) Each  $u_j:\R\to \R, j=1,2,$ is real analytic
and satisfies
\[
\lb{con1}
u'_2(t)\le 0, \qqq \forall \  t\in \R.
\]
(2) There exist nondecreasing functions
$F_j : [0, \iy ) \to [0, \iy ), \ j=1,2, $ such that}
\[
\lb{con2}
 \|u_1(q)\|\le F_1( \|q'\|),  \quad
  \|u'_2(Q)\|\le F_2( \|q\|),  \quad  q\in \mW_1^0.
\]

\bigskip
With Condition (U) in mind, we write $\Delta_q$ instead of
$\Delta_{q,u}$. Our first main theorem is the following.


\begin{theorem}
\lb{TP} The mapping $P:\mW_1^0\ni q \to p=P(q)\in \mH_0$ given by
\er{dP1}  is a real analytic isomorphism between $\mW_1^0$ and
$\mH_0$. In particular, the  operator ${\pa P\/\pa q}$ has a bounded
inverse for each $q\in \mW_1^0$. Moreover, it has the following
properties.

\noindent
(1) The following  inequalities hold:
\[
\begin{aligned}
\lb{eP1}
&\|q'\|\le \|p\| \le \|q'\|+\|q^2\|+
\|q'\|F_1(\|q'\|)+\|q\| F_2(\|q\|)+\|q^2\|^{1\/2}F_2(\|q\|)^{1\/2}.
\end{aligned}
\]
\noindent
(2) The mapping $P(q)-q' : \mW_1^0 \to \mH_0$ is compact.

\noindent
(3) The mapping $\mW_1^{0,odd} \ni q\to p=P(q)$  is a real analytic
isomorphism between $\mW_1^{0,odd}$ and $\mH_0^{even}$.

\noindent
(4) The mapping $\mH_1^{odd} \ni q\to p=P(q)$ is a real analytic
isomorphism between $\mH_1^{odd}$ and $\mH_0^{even}$.

\end{theorem}

\noindent {\bf Remark.} (1)  The mapping $q\to p=q'+q^2+u-c_0 :
\mH_1 \to \mH_0$ was considered  in \cite{K02}. In some cases the
mapping $\mH_0$ into  $\mH_{-1}$ is also useful (see \cite{K03},
\cite{BKK03}).

\noindent (2) For the surface of revolution (\cite{IK13}), we need
to study  the case $u=E\r^{-{4\/d}}$. Here $d+1\ge 2$ is the
dimension of the surface of revolution and $E\ge 0$ is a constant.


\subsection{Second main result - Inverse problems}
We consider the eigenvalue problems for $-\Delta_q$ and $S_p$ on
$(0,1)$  subject to the boundary condition
\begin{equation}
f'(0) - af(0) = 0, \qq f'(1) + bf(1) = 0, \qq
a, b \in  \R\cup\{\infty\}.
\label{S1BC}
\end{equation}
Our second main theorem asserts that the mapping in Theorem \ref{TP}
preserves the boundary conditions and spectral data.

\begin{theorem}
\lb{T2} Let $p=P(q),  \ q\in \mW_1^0$, be defined by \er{dP1}. Then
the operators $S_p$ and $-\D_q$ are unitarily equivalent. In
particular, they
 have the same boundary conditions, eigenvalues and the norming constants.
\end{theorem}

Therefore, the inverse problem for $-\Delta_q$ is solvable if and
only if  so is for $S_p$. Let us consider the  following three cases
separately.


\subsection{Inverse problem for the Dirichlet boundary condition :
$a=b=\iy$.} Denote by $\m_n=\m_n(q),\ n\ge 1$, the eigenvalues of
$-\D_q$ subject to the boundary condition (\ref{S1BC}) for the case
$a =b =\infty$.  It is well-known that all $\m_n$ are simple and
satisfy
$$
\m_n=\m_n^0+c_0+\wt\m_n,\quad {\rm where}\quad
(\wt\m_n)_{1}^{+\iy}\in\ell^2, \qq c_0=\int_0^1(q^2+u)dt,
$$
and $\m_n^0=(\pi n)^2$, $n\ge 1$, denote the unperturbed eigenvalues.
The norming constants are defined by
\[
\label{nc00} \vk_n(q)=\log\left|\r(1)f_n'(1,q)\/f_n'(0,q)\right|,\qquad n\ge 1,
\]
where $f_n$ is the $n$-th eigenfunction satisfying $f_n'(0)\ne 0$
and $f_n'(1)\ne 0$. Applying Theorems \ref{TP} and \ref{T2} and the
result of the inverse problem for $S_p$ in \cite{PT87}, we have the
following theorem.


\begin{theorem}
\label{Tip1}
Let $a=b=\iy$ and let
$\cH=\mW_0^1$ or $\cH=\mH_1$. Then  the mapping
$$
\P : q\mapsto \left((\wt\m_{n}(q))_{n=1}^{\iy}\,;
(\vk_{n}(q))_{n=1}^{\iy}\right)
$$
is a real-analytic isomorphism between $\cH$ and $\cM_1\ts \el2_1$, where
\begin{equation}
\cM_1 = \left\{(h_n)_{n=1}^{\infty}\in \ell^2\, ; \, \mu_1^0 + h_1
 < \mu_2^0 + h_2 < \cdots\right\} \subset \ell^2.
\label{S1DefinecM1}
\end{equation}

In particular, in the anti-symmetric case, i.e., $q$ is odd,
the spectral mapping
\[
\wt\s: \cH^{odd}\to \cM_1,\qqq {\rm given \ by }   \qqq p\to
\wt\s
\]
is a real analytic isomorphism between
$\cH^{odd}$ and  $\cM_1$, where $\cH^{odd}=\mH_0^{odd}$ or
$\cH^{odd}=\mW_0^{1,odd}$.

\end{theorem}


\subsection{Inverse problem for the mixed boundary condition : $a=\iy,\ b\in \R$.}
Let $\l_n=\l_n(q,b), \ n\geq 0$, be the eigenvalues of $-\D_q$
subject to the boundary condition (\ref{S1BC}) for the case $a=\iy,
b\in \R$. We then have
$$
\l_n=\l_n^0+c_0+2b+\wt\l_n(q,b),\quad \mathrm{where} \quad
(\wt\l_n)_{0}^{\iy}\in\ell^2,\qq c_0=\int_0^1(q^2+u)dt,
$$
and $\l_n^0=\pi^2(n+{1\/2})^2$, $n\ge 0$, denote the unperturbed eigenvalues.
The norming constants are defined by
\[
\label{ncb} \c_n(q,b)=\log\left|\r(1)f_n(1,q,b)\/f_n'(0,q,b)\right|,
\qquad n\ge 0,
\]
where $f_n$ is the $n$-th  eigenfunction satisfying $f_n'(0,q,b)\ne
0$  and $f_n(1,q,b)\ne 0$.
 A simple calculation gives
$$
\c_n^0=\c_n(0,0)=-\log \pi(n\!+\!{\textstyle{1\/2}}),\qquad {\rm where}\qquad
\sqrt{\l_n^0}=\pi(n\!+\!{\textstyle{1\/2}}).
$$

Denote by $\vp(x)=\vp(x,\l,q,b)$ and $\x(x)=\x(x,\l,q,b)$  the
fundamental solutions of the equation
\[
\lb{eqq}
-{1\/ \r^2}(\r^2f')'+u f=\l f
\]
 such that
\[
\vp(0)=0, \qqq \vp'(0)=1, \qq  \x(1)=-1,\qq \x'(1)=b.
\]
We define the Wronskian $w(\l,q,b)$ for the equation \er{eqq} by
\[
w=\{\vp, \x  \}_\r=\r(1)\big(\vp'(1,\l,q,b)+b \vp(1,\l,q,b)\big),
\]
where $\{f, g\}_\r=f (\r g')-(\r f') g$.  Note that $\l_n(q,b)$  is
a simple root of $w(\l,q,b)$, and is given by
\[
\label{adamA} w(\l,q,b)=
\cos\sqrt\l\cdot\prod_{n=0}^{+\iy}{\l-\l_n(q,b)\/\l-\l_n^0}\,,\qquad
\l\in\C.
\]
The inverse problem for $S_p$ with $a=\iy, b\in \R$ was solved in
\cite{KC09}. Therefore, applying Theorems \ref{TP},  \ref{T2} and
the result  of the inverse problem for $S_p$ in \cite{KC09}, we have
the following theorem.

\begin{theorem}
\lb{Tip2}
(1) For each $b\in \R$ the mapping
$$
\P:q\mapsto
\left((\wt\l_n(q,b))_{n=1}^{\iy}\,;(\c_{n-1}(q,b)-\c_{n-1}^0)_{n=1}^{\iy}\right)
$$
is a real-analytic isomorphism between $\mW_1^0$ and $\L_1\ts\ell^2_1$, where
$\L_1$ is given by
\[
\lb{L1}
\L_1=\left\{(h_n)_{n=1}^{\iy}\in\ell^2:\l_1^0\!+\! h_1\!<\!\l_{2}^0\!+\! h_{2}\!
<\dots \right\} \ss\ell^2.
\]
(2) For each $(q,b)\in \mW_0^1\ts\R$ the following identity holds true:
\[
\label{IdentityB}
b=\sum_{n=0}^{\iy} \lt(2-{e^{\c_n(q,b)}\/|{\pa w\/\pa \l}(\l_n,q,b)|}\rt),
\]
where the series converges uniformly.
\end{theorem}


\subsection{Inverse problem for the generic boundary condition : $a,b\in \R$.}
Let $\m_n=\m_n(q,a,b), \ n\geq 0,$ be the eigenvalues of $-\Delta_q$
subject to the boundary condition (\ref{S1BC}) for the case $a, b\in
\R$.   Then we have
$$
\m_n=\m_n^0+c_0+2(a+b)+\wt\m_n(q,a,b),\quad \mathrm{where} \quad
(\wt\m_n)_{1}^{\iy}\in\ell^2,\qq c_0=\int_0^1(q^2+u)dt,
$$
and $\m_n^0=(\pi n)^2$, $n\ge 1$, denote the unperturbed eigenvalues.
We introduce the norming constants
\[
\label{ncab} \f_n(q,a,b)=\log\left|\r(1)f_n(1,a,q,b)\/f_n(0,q,a,b)\right|,
\qquad n\ge 0,
\]
where $f_n$ is the $n$-th normalized eigenfunction satisfying
$f_n(1,a,q,b)\ne 0$ and $f_n(0,q,a,b)\ne 0$. The inverse problem for
$S_p$ with  generic  boundary condition was solved in \cite{KC09}.
Therefore, by Theorems \ref{TP},  \ref{T2} and the results in
\cite{KC09}, we have the following theorem.


\begin{theorem}
\label{Tip3} For any $a,b\in\R$, the mapping
$$
\P_{a,b}:q\mapsto \left((\wt\m_{n}(q,a,b))_{n=1}^{\iy}\,;
(\f_{n}(q,a,b))_{n=1}^{\iy}\right)
$$
is a real-analytic isomorphism between $\mW_1^0$ and $\cM_1\ts \el2_1$, where
$\cM_1$ is given by \er{S1DefinecM1}.
\end{theorem}

When $a,b$ are not fixed, the inverse problem for $S_p$ was solved in \cite{IT83}. Therefore,  Theorems \ref{TP}, \ref{T2} and  \cite{IT83} prove the following theorem.

\begin{theorem}
\label{Tip4}
(1) The mapping
$$
\P : (q,a,b)\mapsto \left((\wt\m_{n-1}(q,a,b))_{n=1}^{\iy}\,;
(\f_{n-1}(q,a,b))_{n=1}^{\iy}\right)
$$
is a real-analytic isomorphism between $\R^2\ts \mW_1^0$ and $\cM_0\ts \el2_1$,
where
\[
\lb{dM0}
\cM_0=\left\{(h_{n-1})_{n=1}^{\iy}\in\el2:\m^0_0+h_0\!<\!\m_1^0+
h_1\!<\!\dots\right\}\ss\el2.
\]
(2) For each $(q;a;b)\in \mW_0^1\ts\R^2$ the following identities hold true:
\begin{equation}
\label{Iab} b=-1+\sum_{n=0}^{+\iy} \left(2-{e^{\f_n}\/|{\pa w_2\/\pa
\l}(\m_n)|}\right), \qquad
a=-1+\sum_{n=0}^{+\iy}\left(2-{e^{-\f_n}\/|{\pa w_2\/\pa
\l}(\m_n)|}\right),
\end{equation}
where  $\displaystyle{w_2(\l)\!=-\sqrt\l\sin\sqrt\l\,\cdot
{\prod}_{0}^{\iy} {{\l-\m_n\/\l-\m_n^0}},\qqq \m_n=\m_n(q,a,b),
\f_n=\f_n(q,a,b).}$

\end{theorem}

\subsection{Plan of the paper}
Theorem \ref{TP} is proved in Section 2. Proof of Theorems
\ref{T2}-- \ref{Tip3} will be given in Section 3. In Section 4, we
shall review the results for the 1-dimensional inverse problems for
the Schr{\"o}dinger operator $S_p y=-y''+py$ on the finite interval.

\section {Non-linear mappings}
\setcounter{equation}{0}

\subsection {Unitary transformations.}
 We define a unitary transformation $\mU$ by
\begin{equation}
 \mU:
L^2([0,1],\r dx)\ni f\to \rho f \in L^2([0,1],dx),\qqq \r=e^{\int_0^xq(t)dt},
\nonumber
\end{equation}
which transforms the operator $-\D_q$ into the Schr\"odinger
operator  $S_p$ as follows
\[
\lb{5}
 \mU (-\D_q) \mU^{-1}=
 -\r^{-1}\pa_x \r^2\pa_x r^{-1}+u=\cD^*\cD+u=S_p+c_0.
\]
The operator $S_p$ acts in $ L^2([0,1],dx)$   and is given by
\[
\begin{aligned}
\lb{6}
S_p =-\partial_x^2+p,\qqq c_0=\int_0^1(q^2+u)dx,\\
 \qqq p=q'+q^2+u-c_0,\qqq  \qqq \r=e^{Q}.
\end{aligned}
\]
Here, using  $\r=e^{\int_0^xq(t)dt}$ and letting $A^{\ast}$ denote
the formal adjoint of $A$, we have computed
\[
\begin{aligned}
\lb{a7}
\cD =\r \pa_x \r^{-1}= \pa_x -q,\qqq \cD^* =
\big(\r \pa_x \r^{-1}\big)^*=-\pa_x-q,\\
\cD^*\cD=-(\pa_x +q)(\pa_x -q)=-\pa_x^2+q'+q^2.
\end{aligned}
\]

In this section we consider the mapping $P: \mW_1^0  \to \mH_0$ defined by
\[
\lb{dP3}
\begin{aligned}
\ca  &P(q)=q'+q^2+u-c_0, \qqq  \\
&u=u_1(q)+u_2(Q),\qqq Q=\int_0^x q(t)dt,\  \\
& c_0=\int_0^1 (q^2+u)dt\ac,
\end{aligned}
\]
where $u$ satisfies Condition U.


\subsection { Estimates.} We first derive preliminary estimates on
 the non-linear mapping. Recall that
$$
\|f\|^2 = (f,f), \quad (f,g)=\int_0^1f\ol gdx.
$$
It is easy to show the following inequalities for $q\in \mW_1^0$:
\[
\lb{qdq}
 \|q\|_{L^\iy(0,1)} \le  \|q'\|,\qqq \qqq \|q\| \le  \|q'\|,
\]
\[
\begin{aligned}
\lb{Qq}
\|Q\|_{L^\iy(0,1)} \le  \|q\|.
 \end{aligned}
\]

We need the estimates of $P(q)$ from above and below.

\begin{lemma}\label{Lemma 2.1}
\lb{Tpqw} Let $p=P(q), \, q\in \mW_1^0$ be given by \er{dP3} and
let $h=q^2+u-c_0$. Then the following estimates hold true:
\[
\lb{Pe1} \|q'\|^2 + \|h\|^2 \le  \|p\|^2=\|q'\|^2+\|h\|^2-(q^2,u_2'(Q)),
\]
\[
\lb{Pu2}
\begin{aligned}
\|u_2(Q)- u_2(0)\|\le \|Q\|_{L^\iy(0,1)} F_2(\|q\|),\\
\|u_1(q)- u_1(0)\|\le \|q\|_{L^\iy(0,1)} F_1(\|q'\|),
\end{aligned}
\]
\[
\lb{Pu}
\begin{aligned}
\|u-\int_0^1udx\|\le \|u-u(0)\|\le F(\|q'\|,\|q\|),\\
\|h\|\le  \|q^2\|+\|u-u(0)\| \le \|q^2\|+  F(\|q'\|,\|q\|)
\end{aligned}
\]
where $F(\|q'\|,\|q\|)=\|q'\|F_1(\|q'\|)+\|q\| F_2(\|q\|)$
\[
\lb{Pe3} \|p\| \le \|q'\|+\|h\|+\|q^2\|^{1\/2}F_2(\|q\|)^{1\/2}.
\]
In particular, if $u_1 = u_2 =0$ and  $c_0=\|q\|^2$, we have
\[
\lb{Pe4} \|q'\|^2\le \|p\|^2=\|q'\|^2+\|q^2-c_0\|^2
=\|q'\|^2+\|q^2\|^2-c_0^2 \le  \|q'\|^2+\|q^2\|^2.
\]
\end{lemma}
\no {\bf Proof.} Let $h=q^2+u-c_0$, where $u=u_1(q)+u_2(Q)$. We have
\[
\lb{pq1}
\begin{aligned}
\|p\|^2=\|q'\|^2+\|h\|^2+2(q',q^2+u-c_0)=\|q'\|^2+\|h\|^2+2(q',u),\\
 \end{aligned}
\]
since
$$
(q',1)=0,\qqq (q',q^2)=0.
$$
Letting $U(x) = \int_0^xu_1(t)dt$, we have by integration by parts
and  the Condition $U$
\begin{equation}
\begin{split}
(q',u) & =\int_0^1q'u_1(q)dx+\int_0^1q'u_2(Q)dx\\
 &=U_1(q(x))\Big|_0^1
-\int_0^1q^2u_2'(Q)dx=-(q^2, u_2'(Q))\ge 0,
\end{split}
\nonumber
\end{equation}
which together with \er{pq1} yields \er{Pe1}.

We have the identity
$$
u_2(Q(x))- u_2(0)= \int_0^1Q(x)u_2'(tQ(x))dt.
$$
The inequality \er{Pu2} then  follows from the computation
$$
\begin{aligned}
&\|u_2(Q)- u_2(0)\|^2=\int_0^1(u_2(Q(x))- u_2(0))^2dx=
\int_0^1Q^2(x)dx\rt|\int_0^1u_2'(tQ(x))dt\rt|^2\\
&\le \int_0^1Q^2(x)dx\int_0^1u_2'(tQ(x))^2dt\le \|Q\|_{L^\iy(0,1)}^2
\int_0^1dt\|u_2'(tQ)\|^2\\
&\le \|Q\|_{L^\iy(0,1)}^2 \int_0^1 F_2^2(t\|q\|)dt\le
\|Q\|_{L^\iy(0,1)}^2 F_2^2(\|q\|),
\end{aligned}
$$
where we have used (\ref{con2}) and the monotonicity of $F_2$.

The identity
$$
u_1(q(x))- u_1(0)= \int_0^1q(x)u_1'(tq(x))dt.
$$
and similar arguments ((\ref{con2}) and the monotonicity of $F_2$) yield
$$
\begin{aligned}
&\|u_1(q)-u_1(0)\|^2=\int_0^1(u_1(q(x))- u_1(0))^2dx=
\int_0^1q^2(x)dx\rt|\int_0^1u_1'(tq(x))dt\rt|^2\\
&\le \int_0^1q^2(x)dx\int_0^1u_1'(tq(x))^2dt\le \|q\|_{L^\iy(0,1)}^2
\int_0^1dt\|u_1'(tq)\|^2\\
&\le \|q\|_{L^\iy(0,1)}^2 \int_0^1 F_1^2(t\|q'\|)dt\le
\|q\|_{L^\iy(0,1)}^2 F_1^2(\|q'\|),
\end{aligned}
$$
which gives \er{Pu2}.

Using \er{Pu2} and (\ref{Qq}), we obtain
$$
\|u-u(0)\|\le \|u_1(q)-u_1(0)\|+\|u_2(Q)-u_2(0)\|\le
\|q'\|F_1(\|q'\|)+\|q\|F_2(\|q\|).
$$
This estimate and a simple bound for $v\in L^2(0,1)$
\[
\|v-\int_0^1vdx\|\le \|v\|
\]
give
$$
\begin{aligned}
\|u-\int_0^1udx\|=
\|u-u(0)-\int_0^1(u-u(0))dx\|\le\|u-u(0)\|\le F(\|q'\|,\|q\|),
\end{aligned}
$$
and
$$
\begin{aligned}
\|h\|=\|(q^2-\|q\|^2)+u-u(0)-\int_0^1(u-u(0))dx\|\\
\le
\|(q^2-\|q\|^2\|+\|u-u(0)\|\le \|q^2\|+
F(\|q'\|,\|q\|)
\end{aligned}
$$
which yields \er{Pu}.

Using (\ref{con2}), \er{Pu} and \er{Pe1}, we obtain \er{Pe3}.
If $u_1 = u_2=0$, the inequality \er{Pe4} follows from \er{Pe1}-\er{Pe3}.
\BBox


\subsection { Analyticity and invertibility.} We show that that mapping
$P: \mW_1^0\to \mH_0$ is real analytic.


\begin{lemma}\lb{TApq}
(1)  The map $P: \mW_1^0\to \mH_0$  is
real analytic  and its gradient is given by
\[
\lb{ap1}
 {d P(q)\/d q} f=f'+2q f+u_1'(q)f+u_2'(Q) J f-{d c_0(q)\/d q}f,
 \qqq \ \forall \ q,f\in \mW_1^0,
\]
\[
\lb{ap2} {d c_0\/d q}f=\int_0^1\Big(2qf +u_1'(q)f+u_2'(Q)J  f \Big)dx,
\]
where $Jf=\int_0^xfdx$.

\no
(2) The operator ${d P(q)\/d q} $ is invertible for all
$q\in \mW_1^0$.
\end{lemma}

\no {\bf Proof}. The proof of (1) is standard. See \cite{PT87}.

We show (1) by a contradiction.
 Due to \er{ap1}, the linear operator ${d P(q)\/d q}: \mW_1^0\to
\mH_0$ is a sum of ${d\/dx}$ and a
compact operator for all $q\in \mW_1^0$. Thus it is
a Fredholm operator. We prove that ${d P(q)\/d q}$
is invertible. Let $q\in \mW_1^0$ and $f\in \mW_1^0$ be a solution of
the equation
\[
\lb{ap3x} {d P(q)\/d q}f=0,
\]
which is rewritten as
\[
\lb{ap3}
f'+2q f+u_1'(q)f+u_2'(Q) J f-{d c_0(q)\/d q}f=0.
\]
Setting $y=J f$, we obtain the equation
for $y$ given by
\[
\lb{ap4}
\begin{aligned}
-y''-2wy'+ V y=C,\qqq C=-\int_0^1\Big(2qf +u_1'(q)f+u_2'(Q)J  f \Big)dx,\\
w=q+{1\/2}u_1'(q),\qqq V=-u_2'(Q)\ge 0,\qqq y'(0)=y'(1)=y(0)=0.
\end{aligned}
\]
We rewrite  \er{ap4} in the following form
\begin{equation}
\lb{ap5}
\begin{split}
& -{1\/r^2}(r^2y')'+Vy=C, \quad r(x)=e^{\int_0^x wdt},\\
& y'(0)=y'(1)=y(0)=0.
\end{split}
\end{equation}

If $C=0$, by multiplying
\er{ap5} by $r^2y$, we have
$$
0=\int_0^1\Big(-(r^2y')'y +r^2 V y^2\Big)dx= \int_0^1\Big((r y')^2 +
r^2V y^2\Big)dx,
$$
which implies $y=0$, hence $f = 0$.

If $C\ne 0$, we can assume without loss of generality that
$C=1$ in \er{ap5}.
The solution of \er{ap5} with $C=1$ has the form
\begin{equation}
y(x)=\int_0^1R(x,t)r^2(t)dt\ge 0,
\label{y(x)=intRxt}
\end{equation}
 where $R(x,t)$ is the Green function for the problem
\begin{equation}
- \frac{1}{r^2}(r^2\psi')' + V\psi = f, \quad 0<x<1, \quad
y(0) =y'(1)=0.
\label{StLiouprob}
\end{equation}
 Note that
 $R(x,t)\ge 0$ for any $x,t\in [0,1]$, since the first eigenvalue
  of the Sturm-Liouville problem
equation
\[
\lb{ap6}-{1\/r^2}(r^2\p')'+V\p=\l \p, \qqq \p(0)=\p'(1)=0,
\]
is positive.
Let $\vp (x,t)$ be the
solution of the equation with parameter $t \in [0,1]$:
$$
-\vp''-2w(x+t)\vp'+V(x+t)\vp =0 ,\qqq \vp (0,t)=0, \qq \ \vp'
(0,t)=1.
$$
A direct calculation shows that $z(x)$ defined by
$$
z(x)=-\int_0^x\vp (x-t,t)dt
$$
satisfies
$$
-z''-2w(x)z' + V(x)z = 1, \quad z(0)=z'(0)=0,
$$
hence $z(x)=y(x)$. Note that the eigenvalues $\lambda$ of the
Sturm-Liouville problem (\ref{ap6}) with $r(x), V(x)$ replaced by
$r(x+t), V(x+t)$ are positive. Then by the well-known comparison
principle, $\vp (x,t)$ is positive for $x,t\in [0,1]$.  Therefore,
$y(x) < 0$ on $[0,1]$. In view of \er{y(x)=intRxt}, we have arrived
at a contradiction. $\BBox$

\subsection {Bijectivity of $P$.}

The proof of  Theorem \ref{TP} is done by the "direct approach"
 in  \cite{KK97} based on nonlinear functional analysis, the
main theorem of which is the following one.

\begin{theorem}
\lb{TA97}
 Let $H, H_1$ be real separable Hilbert spaces equipped with norms
$\|\cdot \|, \|\cdot \|_1$. Suppose that the map $f: H \to H_1$
satisfies the following conditions:

 (i) $f$ is real analytic  and the operator ${df /dq}$ has a bounded
inverse for all $q\in H$,


 (ii) there is a nondecreasing function $\e: [0, \iy ) \to [0, \iy
)$ such that $\e (0)=0$ and $\|q\|\le \e (\|f(q)\|_1)$  for all $q\in
H$,

(iii) there exists a linear isomorphism $f_0:H\to H_1$ such that
the map $f-f_0: H \to H_1$  is compact.

\no Then $f$ is a real analytic isomorphism between $H$ and $H_1$.
\end{theorem}

We now prove the main results of this section.

\medskip

{\bf Proof of Theorem \ref{TP}.} We check all conditions in Theorem
\ref{TA97} for the mapping $P(q), q\in \mW_1^0$. The condition  (i)
is proved in Lemma \ref{TApq}.  The two-sided  estimates \er{eP1}
are proved in Lemma \ref{Tpqw}. Let $q^\n\to q$ weakly in $\mW_1^0$
as $\n\to \iy$. Then we have that $q^\n\to q$ strongly in $\mH_0$ as
$\n\to \iy$, since the imbedding $\mW_1^0\to\mH_0$ is compact. Hence
the mapping $q\to P(q)-q'$ is compact, and the condition (iii) is
satisfied. Therefore, by Theorem \ref{TA97},
 $P:\mW_1^0 \to \mH_0$ is a real analytic isomorphism  between
 the Hilbert spaces $\mW_1^0$ and $\mH_0$.

If $q\in \mW_1^{0,odd}$, then $p=P(q)\in \mH_0^{even}$,  by the
definition of $P$. Repeating the arguments for $P(q)$ when $q\in
\mW_1^{0}$, we deduce that the mapping $P:\mW_1^{0,odd} \to
\mH_0^{even}$ is a real analytic isomorphism  between
 the Hilbert spaces $W_1^{0,odd}$ and $\mH_0^{even}$.

 The proof for the case
 $P:\mH_0^{odd}\to \mH_0^{even}$ is similar.
\BBox


\section {Proof of Theorems \ref{T2}-- \ref{Tip3}}
\setcounter{equation}{0}
 Letting $y=\mU f= \r f$ for $q\in \mW_0^1$, we obtain
\[
\lb{bc1}
\begin{aligned}
& y(0)=f(0),\qqq y'(0)=f'(0),\\
& y(1)=\r(1)f(1),\qqq y'(0)=\r(1)f'(1).
\end{aligned}
\]
Thus the function $f$ satisfies the  boundary condition
\[
\lb{bc2}
f'(0)-af(0)=0,\qquad f'(1)+b f(1)=0,\qquad a,b\in \R\cup \{\iy\},
\]
if and only if $y=\mU f= \r f$   satisfies the same boundary condition
\[
\lb{bc3}
y'(0)-ay(0)=0,\qquad y'(1)+by(1)=0,\qquad a,b\in \R\cup \{\iy\},
\]
since $q\in \mW_1^0$ and $\r'(0)=\r'(1)=0$.

\subsection {Proof of Theorem \ref{T2}}

Let $p=P(q), q\in W_1^0$, be defined by \er{dP1}. Then under the
transformation $y=\mU f=\r f$ the operators $S_p$ and $-\D_q$ are
unitarily equivalent. Moreover, due to \er{bc1}-\er{bc3} the
operators $S_p$ and $-\D_q$ have the same the boundary conditions.
The norming constants are the same by \er{bc1}.

Assume that the mapping $p\to $ "eigenvalues + norming constants for
the operator $S_p$" gives the solution of the inverse problem for
the operator $S_p$. Then, since the mapping $p\to q$ is an analytic
isomorphism, we obtain the solution of the inverse problem for
$-\Delta_q$ from that of $S_p$, and vice versa. \BBox

\medskip
In the proof of Theorems \ref{Tip1}-\ref{Tip3}, we give a more
detailed  explanation.

\subsection{Proof of Theorems \ref{Tip1}-\ref{Tip3}.}
Recall that in \er{5} the operator $-\D_q$ acting in
$L^2([0,1],\r^2dx)$ is shown to be unitarily equivalent to
\[
\lb{a6}
\mU (-\D_q) \mU^{-1} =S_p+c_0,
\]
acting in $L^2([0,1],dx)$.
This representation is more convenient for us.

\medskip
{\bf Proof of Theorem \ref{Tip1}.}
Let $q\in \mW_0^1$ and $a,b\in\R$.
We consider the Sturm-Liouville problem with the Dirichlet  boundary conditions:
$$
-{1\/\r^2}(\r^2f')'+uf=\l f,\qqq f(0)=f(1)=0.
$$
Let $\m_n=\m_n(q), n=1,2,3,...$ be the eigenvalues.
It is well-known that
$$
\m_n=(\pi n)^2+c_0+\wt\m_n(q),\quad \mathrm{where} \quad
(\wt\m_n)_{1}^{\iy}\in\ell^2,\qq c_0=\int_0^1(q^2+u)dt.
$$
Following \cite{PT87}, we introduce the norming constants
\[
\label{ncb3}
\f_n(q)=\log\left|\r(1)f_n'(1,q)\/f_n'(0,q)\right|,\qquad n\ge 1,
\]
where $f_n$ is the $n$-th eigenfunction. Note that $f_n'(0)\ne 0$.
We then have the mapping
$$
\P: q\mapsto \P(q)=
\left((\wt\m_n(q))_{n=1}^{\iy}\,;(\f_n(q))_{n=1}^{\iy}\right)
$$
We let $p=P(q), q\in \mW_0^1$, and apply Theorem \ref{TKC2}.
Consider  the Sturm-Liouville problem with the Dirichlet  boundary
conditions:
$$
S_p y=-y''+p(x)y,\qqq\qqq y(0)=y(1)=0.
$$
Denote by $\s_n=\s_n(p), n\ge 1$ the eigenvalues of $S_p$ and let
$\vk_n(p)$ be the corresponding norming constants given by
\[
\label{NuDef}
\vk_n(p)=\log\lt|{y_n'(1,p)\/y_n'(0,p)}\rt|\,,\qquad n\ge 0.
\]
Recall that  due to \cite{PT87}  (see Theorem \ref{TMOT})
  the mapping
$$
\F:p\mapsto \F(p)=
\left((\wt\s_{n}(p))_{n=1}^{\iy}\,;(\vk_n(p))_{n=1}^{\iy}\right)
$$
is a real-analytic isomorphism between $\mH_0$ and $\cM_1\ts\el2_1$, where
$\cM_1$ is defined by \er{S1DefinecM1}.
Due to Theorem \ref{T2} we obtain the identity
$$
\F(P(q))=\P(q),\qq \forall \ q\in \mW_1^0.
$$
This is a composition of two mappings $\F$ and $P$,
each of which is an analytic isomorphism (see Theorem \ref{TMOT}
and \ref{TP}). Then  the mapping
$$
\P:q\mapsto
\left((\wt\m_n(q))_{n=1}^{\iy}\,;(\f_n(q))_{n=1}^{\iy}\right)
$$
is a real-analytic isomorphism between $\mW_0^1$ and $\cM_1\ts\ell^2_1$.

Consider the  the spectral mapping in the symmetric case
$$
\wt\m: \mW_1^{0,odd}\to \cM_1,\qqq {\rm given \ by }   \qqq p\to
\wt\m.
$$
Recall that  due to \cite{PT87}  (see Theorem \ref{TMOT})
  the mapping
$$
\wt\s : p\mapsto \wt\s(p)=(\wt\s_{n}(p))_{n=1}^{\iy}
$$
is a real-analytic isomorphism between $\mH_0^{even}$ and $\cM_1$, where
$\cM_1$ is defined by \er{S1DefinecM1}.
Again in view of Theorem \ref{T2}, we obtain the identity
$$
\wt \s(P(q))=\wt \m(q),\qq \forall \ q\in \mW_1^{0,odd}.
$$
This is a composition of two mappings $\wt \s$ and $P$,
each of which being analytic isomorphism (see Theorem \ref{TMOT}
and \ref{TP}). Then so  is  the mapping
$$
\wt\m:q\mapsto  (\wt\m_n(q))_{n=1}^{\iy}  : \mW_1^{0,odd}\to \cM_1.
$$

The proof for the case $q\in \mH_1^{odd}$ is similar.
\BBox

\medskip
{\bf Proof of Theorem \ref{Tip2}.} Let $q\in \mW_0^1$ and $b\in\R$.
We consider the Sturm-Liouville problem with the mixed  boundary condition,
$$
-{1\/\r^2}(\r^2f')'+uf=\l f,\qqq f(0)=0,\qquad f'(1)+b f(1)=0.
$$
Let $\l_n=\l_n(q,b), n=0,1,2,...$ be the associated eigenvalues. Then,
$$
\l_n=\pi^2(n\!+\!{\textstyle{1\/2}})^2+c_0+2b+\wt\l_n(q,b),\quad
\mathrm{where} \quad (\wt\l_n)_{1}^{\iy}\in\ell^2,\qq
c_0=\int_0^1(q^2+u)dt.
$$
As in  \cite{KC09}, the norming constants are defined by
\[
\label{ncb4}
\c_n(q,b)=\log\left|\r(1)f_n(1,q,b)\/f_n'(0,q,b)\right|,\qquad n\ge
0,
\]
where $f_n$ is the $n$-th normalized eigenfunction such that $f_n'(0,q,b)>0$.
A simple calculation gives
$$
\c_n^0=\c_n(0,0)=-\log \pi(n\!+\!{\textstyle{1\/2}}),\qquad {\rm where}\qquad
\sqrt{\l_n^0}=\pi(n\!+\!{\textstyle{1\/2}}).
$$
Thus for fixed $b\in \R$  we have a mapping
$$
\P_b:q\mapsto \P_b(q)=
\left((\wt\l_n(q,b))_{n=1}^{\iy}\,;(\c_{n-1}(q,b)-\c_{n-1}^0)_{n=1}^{\iy}\right).
$$

Let $p=P(q), q\in \mW_0^1$, and apply Theorem \ref{TKC2}. Consider
the Sturm-Liouville problem
$$
S_p y=-y''+p(x)y,\qquad y(0)=0,\qquad y'(1)+by(1)=0,\quad b\in \R.
$$
Denote by $\t_n=\t_n(p,b), n\ge 0$ the eigenvalues of $S_p$ and let
$\n_n(p,b)$ be the corresponding norming constants given by
\[
\label{NuDef4}
\n_n(p,b)=\log\lt|{y_n(1,p,b)\/y_n'(0,p,b)}\rt|\,,\qquad n\ge 0.
\]
Recall that  due to  \cite{KC09}  (see Theorem \ref{TKC2}),
for each $b\in\R$  the mapping
$$
\F_b:p\mapsto \F_b(p)=
\left((\wt\t_{n}(p))_{n=1}^{\iy}\,;(\n_{n-1}(p)-\n_{n-1}^0)_{n=1}^{\iy}\right)
$$
is a real-analytic isomorphism between $\mH_0$ and $\L_1\ts\el2_1$, where
$\L_1$ is defined by \er{L1}.
By Theorem \ref{T2}, we have an identity
$$
\F_b(P(q))=\P_b(q),\qq \forall \ q\in \mW_0^1,
$$
 where  $\F_b$ and $\P$ are  analytic isomorphisms (see Theorem \ref{TKC2}
and \ref{TP}). Then for each $b\in \R$ the mapping
$$
\P_b:q\mapsto
\left((\wt\l_n(q,b))_{n=1}^{\iy}\,;(\c_{n-1}(q,b)-\c_{n-1}^0)_{n=1}^{\iy}\right)
$$
is a real-analytic isomorphism between $\mW_0^1$ and $\L_1\ts\ell^2_1$.
\BBox

\medskip

We consider the mapping $\P_b$,  when $b$ is not fixed.

\begin{theorem}
\lb{Tip2b}

(i) The mapping
$$
\P:(q,b)\mapsto
\left((\wt\l_{n-1}(q,b))_{n=1}^{\iy}\,;(\c_{n-1}(q,b)-
\c_{n-1}^0)_{n=1}^{\iy}\right)
$$
is a real-analytic isomorphism between $\mW_0^1\ts \R$ and $\R\ts\L_0\ts\ell^2_1$, where $\L_0$ is given by
\[
\lb{L0}
\L_0=\left\{(h_{n-1})_{n=1}^{\iy}\in\ell^2:\l_0^0\!+\! h_0\!
<\l_1^0\!+\! h_1\!<\!\l_{2}^0\!+\! h_{2}\!
<\dots \right\} \ss\ell^2.
\]
(ii) For each $(q;b)\in \mW_0^1\ts\R$ the following identity holds true:
\[
\label{IdentityBz} b=\sum_{n=0}^{+\iy}
\lt(2-{e^{\c_n(q,b)}\/|\dot{w}(\l_n,q,b)|}\rt).
\]
\end{theorem}

{\bf Proof.} The proof is similar to the case of fixed $b\in \R$,
and is based on Theorems \ref{TP}, \ref{T2} and \ref{TipbR}. \BBox

\medskip

{\bf Proof of Theorem \ref{Tip3}.}
Let $q\in \mW_0^1$ and $a,b\in\R$.
We consider the Sturm-Liouville problem with the generic  boundary condition,
$$
-{1\/\r^2}(\r^2f')'+uf=\l f,\qqq f'(0)-af(0)=0,\qqq f'(1)+b f(1)=0.
$$
Let $\m_n=\m_n(q,a,b), n=0,1,2,...$ be the associated eigenvalues.
We then have
$$
\m_n=(\pi n)^2+c_0+2(a+b)+\wt\m_n(q,b),\quad \mathrm{where} \quad
(\wt\m_n)_{1}^{\iy}\in\ell^2,\qq c_0=\int_0^1(q^2+u)dt.
$$
As in  \cite{KC09},  the norming constants are defined by
\[
\label{ncb5}
\f_n(q,a,b)=\log\left|\r(1)f_n(1,q,a,b)\/f_n(0,q,a,b)\right|,\qquad
n\ge 0,
\]
where $f_n$ is the $n$-th normalized eigenfunction such that $f_n'(0)>0$.
Thus for fixed $a,b\in \R$  we have the mapping
$$
\P_{a,b}:q\mapsto \P_{a,b}(q)=
\left((\wt\m_n(q,a,b))_{n=1}^{\iy}\,;(\f_n(q,a,b))_{n=1}^{\iy}\right)
$$

As above, we consider the Sturm-Liouville problem
$$
S_p y=-y''+p(x)y,\qquad y'(0)-ay(0)=0,\qquad y'(1)+by(1)=0,\quad a,b\in \R.
$$
Denote by $\s_n=\s_n(p), n\ge 0$ the eigenvalues of $S_p$ and  let
$\vk_n(p)$ be the corresponding norming constants given by
\[
\label{NuDef5}
\vk_n(p)=\log\lt|{y_n(1,p,b)\/y_n'(0,p,b)}\rt|\,,\qquad n\ge 0.
\]
By \cite{KC09}  (see Theorem \ref{TabKC09})
for each $a,b\in\R$  the mapping
$$
\F_{a,b}:p\mapsto \F_{a,b}(p)=
\left((\wt\s_{n}(p))_{n=1}^{\iy}\,;(\vk_n(p))_{n=1}^{\iy}\right)
$$
is a real-analytic isomorphism between $\mH_0$ and $\cM_1\ts\el2_1$, where
$\cM_1$ is defined by (\ref{S1DefinecM1}).
Theorem \ref{T2} gives the following identity
$$
\F_{a,b}(P(q))=\P_{a,b}(q),\qq \forall \ q\in \mW_0^1.
$$
By the same arguments as above, using  Theorem \ref{TabKC09}
and \ref{TP},  for each $a,b\in \R$ the mapping
$$
\P_{a,b}:q\mapsto
\left((\wt\m_n(q,a,b))_{n=1}^{\iy}\,;(\f_n(q,a,b))_{n=1}^{\iy}\right)
$$
is a real-analytic isomorphism between $\mW_0^1$ and $\cM_1\ts\ell^2_1$.
\BBox

\section {Review of 1-dimensional inverse problems}
\setcounter{equation}{0}

\subsection {\bf Results of Gel'fand-Levitan-Marchenko-Ostrovski-Trubowitz}
 Let us briefly review the inverse spectral
theory for Strum-Liouville operators on a finite interval. We recall only some
important steps mostly focusing on the {\it characterization} problem, i.e.,
the complete description of spectral data that correspond to some fixed class
of potentials. More information about different approaches to inverse spectral
 problems can be found in the monographs \cite{M77}, \cite{L84},
 \cite{PT87} and references therein.

The inverse spectral theory goes back to the seminal \mbox{paper
\cite{Bo46}} (see also a simpler proof \cite{Le49}). Borg showed
that spectra of two Sturm-Liouville problems $-y''+p(x)y=\l y, x\in
[0,1]$, with the same boundary conditions at $x=1$ but different
boundary  conditions at $x=0$, determine the potential $p(x)$ and
the boundary conditions uniquely. The next step was done by
Marchenko \cite{M50}. He  proved that the so-called spectral
function (or, equivalently, the Weyl-Titchmarsh function) determines
 the potential uniquely. Note that the spectral function is
 piecewise-linear outside the spectrum $\{\l_n\}_{n=1}^{+\iy}$
and its jump at $\l_n$ is equal to the so-called {\it normalizing constant}
${1/\a_n(p)}$
given \er{ncM}.
 At the same time, a different approach to this problem was
developed by Krein \cite{Kr51}, \cite{Kr53}, \cite{Kr54}.

An important result was obtained by Gel'fand and Levitan \cite{GL51}.
They gave an effective method to reconstruct the potential $p(x)$
 from its spectral function. More precisely, they derived an integral
  equation and expressed $p(x)$ explicitly in terms of the solution of
this equation. At that time, there was some gap between necessary and
 sufficient conditions
for the spectral functions corresponding to fixed classes of $p(x)$.

Some characterization of spectral data for $q$ such that
$p^{(m)}\in L^1(0,1)$ was derived by Levitan and Gasymov \cite{LG64}
for all $m=0,1,2,..$. Also, they formulated the solution of
the characterization problem in the case $p''\in L^2(0,1)$ (without proof).
Marchenko and Ostrovski \cite{MO75} obtained a sharpening of this result.
Namely, for all $m=0,1,2,..$ they gave the complete
solution of the inverse problem in terms of two spectra,
if $p^{(m)}\in L^2(0,1)$.

So, the inverse problem in the case of Dirichlet boundary conditions
\er{H1} a sharp characterization of  all spectral data
$(\s_n(p),\a_n(p))_{n=1}^{+\iy}$\,) that correspond to potentials
$p\in\mH_0$ is available due to \cite{MO75}. Namely, the necessary
and sufficient conditions are
\begin{equation}
\label{xMO75}
\begin{array}{cl}
\s_1<\s_2<\s_3<...,\qquad & (\s_n-\pi^2n^2)_{n=1}^{+\iy}\in\ell^2 \cr
\mathrm{and}\vphantom{|^\big|} &
(2\pi^2n^2\a_n(p)-1)_{n=1}^{+\iy}\in\ell_1^2.
\end{array}
\end{equation}

Trubowitz and co-authors (Isaacson \cite{IT83}, Isaacson-McKean
\cite{IMT84},  Dahlberg \cite{DT84}, P\"oschel \cite{PT87})
suggested another approach  to the inverse Sturm-Liouville problem
on the finite interval with separated boundary conditions. It is
based on the analytic properties of the mapping
$\mathrm{\{potentials\}\to\{spectral\ data\}}$ and the explicit
transforms corresponding  to the change of only  a {\it finite}
number of spectral parameters
\mbox{$(\s_n(p),\f_n(p))_{n=1}^{\iy}$}. Their norming constants
$\f_n(p)$ (defined in \er{f4}) differ slightly from  Marchenko's
normalizing constants $\a_n(p)$, see  \er{ncM}, but the
characterizations are equivalent (see Appendix in \cite{CK09}).
Also, this approach was applied to other scalar inverse problems
with purely discrete spectrum:

1) Schr{\"o}dinger operators on the circle (periodic boundary
conditions): Garnett-Trubowitz \cite{GT84}, \cite{GT87},
Kargaev-Korotyaev \cite{KK97} for even potentials and Korotyaev
\cite{K97}, \cite{K99} for general potentials.

2) Impedance equations  on the unit interval Coleman-McLaughlin \cite{CM93},
and an impedance equation  on the  the circle \cite{K00}.

3) Singular Sturm-Liouville operators on $[0,1]$: \cite{GR88}, \cite{Ca97},
\cite{Se07}.

4) Perturbed harmonic oscillators: \cite{MT81}, \cite{CKK04}, \cite{CK07}.

5)  Vector-valued Sturm-Liouville operators
on the unit interval with the Dirichlet boundary conditions \cite{CK09}.

Note that there are results, based on other methods, see e.g.,
\cite{FY02}, \cite{AHM05} for the operator $-\D_q$ and \cite{Ho05}
for the operator  $S_p$ and references therein. Since we use the
analytic method of Trubowitz, we mention only the papers using this
approach.

Thus,  the inverse spectral theory for the Sturm-Liouville operators
is bow well understood. We summarize below the characterization
theorems for the case of the Schr{\"o}dinger operator $- {d^2\/dx^2}
+ p(x)$ on the finite interval.

\subsection{Inverse problem for Dirichlet boundary conditions : $a=b=\iy$.}
We consider the Sturm-Liouville problem with the Dirichlet boundary conditions:
\[
\lb{H1}
\begin{aligned}
S_p y=-y''(x)+p(x)y(x),\qqq y(0)=y(1)=0,
\end{aligned}
\]
where a potential $p$ belongs to the real space $\mH_0$. For further
information see \cite{PT87}. The spectrum of  this operator $S_p$ is
discrete and consists of simple eigenvalues $\s_n, n=1,2,....$,
which satisfy
\[
\lb{SL1}
\begin{aligned}
\s_1<\s_2<\s_3<....< \s_n<\s_{n+1}<.... \\
\s_n=(\pi n)^2+\wt\s_n, \qqq n\ge 1,  \qqq
\wt\s=(\wt\s_n)_1^\iy\in \ell^2.
\end{aligned}
\]
Let $y_{n}(x)=y_{n}(x,p)$ be the corresponding real eigenfunctions, which satisfy  $y_n(0)=y_n(1)=0 $ with the condition $y_n'(0)=1$. Introduce the norming constants $\f_n$ by
\[
\lb{f4}
\f_n=\log |y_n'(1)|,\qqq n\ge 1.
\]
  It is well known that
\[
\lb{f5}
 \f=(\f_n)_1^\iy\in \ell_{1}^2.
\]
The monotonicity \er{SL1} gives that, even if $p$ varies over  the whole space
$\mH_0$,  the image of $\left\{(\wt\s_n)_{1}^{\iy}\right\}$ does
not coincide with  $\ell^2$.
Note that Marchenko \cite{M50} solved the inverse problem in terms of eigenvalues
and {\it normalizing constant} $\a_n(p)$ given
\[
\lb{ncM}
\a_n(p)=\int_0^1y_n^2(x,p)dx.
\]

The results of Marchenko-Ostrovski-Trubowitz are formulated in  the
following theorem.

\begin{theorem}
\lb{TMOT}
The mapping $\F: \mH_0\to \cM_1\ts \ell_1^2$ given by
\[
\lb{pmap}
p\to \F=(\wt\s, \f)
\]
is a real  real analytic isomorphism between the Hilbert space $\mH_0$ and
the set $\cM_1\ts \ell_1^2$.

In particular, in the case of even potentials $p$,
the spectral mapping
\[
\wt\s: \mH_0^{even} \to \cM_1,\qqq {\rm given \ by }   \qqq p\to
\wt\s
\]
is a real  real analytic isomorphism between the Hilbert space
$\mH_0^{even}$ and  $\cM_1$.
\end{theorem}

\subsection{Inverse problem for the  mixed boundary conditions: $a=\iy, b\in \R$.}
Below we describe the results of  \cite{KC09} for the case
$a\!=\!\iy$, $b\!\in\!\R$. As far as we know, this is the most
detailed characterization of spectral data in this case available in
the literature. Let
$$
S_p y=-y''+p(x)y,\qquad y(0)=0,\qquad y'(1)+by(1)=0,\quad b\in \R.
$$
Denote by $\t_n=\t_n(p,b), n\ge 0$ the eigenvalues of $S_p$.  It is
well  known that all $\t_n$ are simple and
$$
\t_n(p,b)=\pi^2(n\!+\!{1\/2})^2+2b+\wt\t_n(p,b),\quad {\rm where}\quad
(\wt\t_n(p,b))_{n=1}^{+\iy}\in\ell^2.
$$
Here  $\pi^2(n\!+\!{1\/2})^2, n\ge 0$, denote the unperturbed eigenvalues.

Let $Y_1(x)=Y_1(x,\l,p,b)$, $Y_2(x)=Y_2(x,\l,p,b)$ be the solutions
of  $-y''+p(x)y=\l y$ such that
$$
Y_1(0)=0,\ \ Y_1'(0)=1\quad {\rm and}\quad Y_2(1)=-1,\ \ Y_2'(1)=b.
$$
Here and below $(\,')={\pa\/\pa x}$ and $(\dot{{\,}})={\pa\/\pa
\l}$.  Note that $\t_n(p,b)$ are the roots of the Wronskian
$$
w_1(\l)=w_1(\l,p,b)\equiv\{Y_1,Y_2\}(\l,p,b)\equiv
Y_1'(1,\l,p,b)+bY_1(1,\l,p,b),\qquad\l\in\C,
$$
where $\{Y_2,Y_1\}=Y_2Y_1'-Y_2'Y_1$. The Hadamard Factorization Theorem implies
\[
\label{adam} w_1(\l,p,b)\equiv
\cos\sqrt\l\cdot\prod_{n=0}^{+\iy}{\l-\t_n(p,b)\/\l-\t_n^0}\,,\qquad \l\in\C.
\]
Let $y_n(x)=y_n(x,p,b)$ be the $n$-th eigenfunction of $S_p$. Note that
$y'_n(0)\!\ne \!0$. We introduce the norming constants as
\[
\label{NuDef42}
\n_n(p,b)=\log\lt|{y_n(1,p,b)\/y_n'(0,p,b)}\rt|\,,\qquad n\ge 0.
\]
A simple calculation gives
$$
\n_n^0=\nu_n(0,0)=-\log \pi(n\!+\!{\textstyle{1\/2}}),\qquad {\rm where}\qquad
\sqrt{\s_n^0}=\pi(n\!+\!{\textstyle{1\/2}}).
$$

Let the boundary parameter $b\in\R$ be fixed (e.g., $b\!=\!0$
corresponds to the boundary conditions $y(0)\!=\!0$, $y'(1)\!=\!0$).
In this case spectral data $(\t_n)_{n=0}^{\iy}$,
$(\n_n)_{n=0}^{\iy}$ are not independent, since they satisfy the
nonlinear equation (\ref{IdentityBx}). It turns out that the first
eigenvalue $\t_0(p,b)$ can be uniquely reconstructed from the other
spectral data $(\t_n)_{n=1}^{\iy}$ and $(\n_n)_{n=0}^{\iy}$.

We formulate the result of \cite{KC09}.

\begin{theorem}
\lb{TKC2}
 Let $b\in\R$ be fixed. Then  the mapping
$$
\F_b:q\mapsto
\left((\wt\t_{n}(p))_{n=1}^{\iy}\,;(\n_{n-1}(p)-\n_{n-1}^0)_{n=1}^{\iy}\right)
$$
is a real-analytic isomorphism between $\mH_0$ and $\L_1\ts\el2_1$, where
$\L_1$ is defined by (\ref{L1}).
 For each $(p;b)\in \mH_0\ts\R$ the following identity holds true:
\[
\label{IdentityBx}  b=\sum_{n=0}^{\iy}
\lt(2-{e^{\n_n(p,b)}\/|\dot{w}(\t_n,p,b)|}\rt).
\]

\end{theorem}

It is possible to ``exclude'' from the spectral data not the first
eigenvalue $\l_0$ but an arbitrary norming constant $\m_m$ (thus,
obtaining the parametrization of iso-spectral manifolds). For any
$\{\l_n^*\}_{n=0}^{\iy}$ such that $\l_n^*=\l_n^0+c^*+\m_n^*$, where
$\left(c^*\,;\{\m_n^*\}_{n=0}^{\iy}\right)\in\R\ts\cM$, denote
$$
w^*(\l):=\cos\sqrt\l\cdot\prod_{n=0}^{\iy}{\l\!-\!\l_n^*\/\l\!-\!\l_n^0},\qquad \l\in\C.
$$
Recall that $w(\l,q,b)\equiv w^*(\l)$ for all
$q\!\in\!\Iso_b\left[\{\l_n^*\}_{n=0}^{\iy}\right]$.
\begin{corollary}\label{CorFixedBNC}
Let $b\in\R$ and
$\left(c^*;\{\m_n^*\}_{n=0}^{\iy}\right)\in\R\ts\cM$.  Then for each
$m\ge 0$ the mapping
$$
q\mapsto \{\n_n(q,b)-\n_n^0\}_{n=0,n\ne m}^{\iy}
$$
is a real-analytic isomorphism between the isospectral set
$\Iso_b\left[\{\l_n^*\}_{n=0}^{\iy}\right]$ (which is a real-analytic
 submanifold of
$L^2(0,1)$) and the open set
\[
\label{NbmDef} \cN^b_m= \lt\{\{\n_n-\n_n^0\}_{n=0,n\ne m}^{\iy}\in\el2_1:
\sum_{n=0,n\ne
m}^{\iy} \lt(2-{e^{\n_n}\/|\dot{w}^*(\l_n^*)|}\rt)>b-2\rt\}\ss\el2_1.
\]
\end{corollary}

We consider the case when $b$ is not fixed.

\begin{theorem}
\label{TipbR}
 The mapping
$$
\F:(p;b)\mapsto
\left((\wt\t_{n-1}(q,b))_{n=1}^{\iy}\,;(\n_{n-1}(p,b)-\n_n^0)_{n=1}^{\iy}\right)
$$
is a real-analytic isomorphism between $\mH_0\ts\R$ and $\cM_0\ts\el2_1$,
where $\cM_0$ is given by \er{dM0}.
\end{theorem}

\subsection{Inverse problem for generic boundary conditions: ${\bf a,b\in \R}$.}
Let $\s_n=\s_n(p,a,b)$, $n=0,1,2,...$ be the eigenvalues of  the
Sturm-Liouville problem
$$
-y''+p(x)y=\l y,\qquad y'(0)-ay(0)=0,\qquad y'(1)+b y(1)=0,
$$
where $p\in \mH_0$ and $a,b\in\R$. It is well known that
$$
\s_n=\pi^2n^2+2(a+b)+\wt\s_n(p,a,b),\quad \mathrm{where} \quad
(\wt\s_n)_{1}^{\iy}\in\ell^2,
$$
 where $\pi^2n^2$, $n\ge 0$, denote the unperturbed eigenvalues. Note that
  $\s_n$ are the (simple) roots of the Wronskian
$$
w(\l,p,a,b)\equiv (\vt'\!+\!a\vp'+b(\vt\!+\!a\vp))(1,\l,p),\qquad \l\in\C.
$$
Following \cite{IT83}, we introduce the norming constants
\[
\label{SnDef}
\vk_n(p,a,b)=\log\left|y_n(1,p,a,b)\/y_n(0,p,a,b)\right|,\qquad n\ge 0,
\]
where $y_n$ is the $n$-th normalized eigenfunction such that $y_n(0)>0$.

The  case for fixed  $a,b$ is important.
We recall the results of  \cite{KC09}.

\begin{theorem}
\label{TabKC09} For any $a,b\in\R$ the mapping
$$
\P_{a,b}:p\mapsto \left((\wt\s_{n}(p,a,b))_{n=1}^{\iy}\,;
(\vk_{n}(p,a,b))_{n=1}^{+\iy}\right)
$$
is a real-analytic isomorphism between $\mH_0$ and $\cM_1\ts \el2_1$,
 where $\cM_1$ is given by (\ref{S1DefinecM1}).
\end{theorem}

For the case when $a,b$ are not fixed, we refer the result
of \cite{IT83}.

\begin{theorem}
\label{TipabR}
The mapping
$$
\P:(p;a;b)\mapsto
\left( (\wt\s_{n-1}(p,a,b))_{n=1}^{\iy}\,;(\vk_{n-1}(p,a,b))_{n=1}^{\iy}\right)
$$
is a real-analytic isomorphism between $\mH_0\ts\R^2$ and $\cM_0\ts\el2_1$, where
$$
\cM_0=\left\{(h_{n-1})_{n=1}^{\iy}\in\el2:\m^0_0+h_0\!<\!\m_1^0+
h_1\!<\!\dots\right\}\ss\el2.
$$
\end{theorem}

\bigskip

\setlength{\itemsep}{-\parskip} \footnotesize \no  {\bf
Acknowledgments.} Various parts of this paper were written during
Evgeny Korotyaev's stay in the Mathematical Institute of University
of Tsukuba, Japan and  Mittag-Leffler Institute, Sweden. He is
grateful to the institutes for the hospitality. His study was
supported by the Ministry of education and science of Russian
Federation, project   07.09.2012  No  8501 No
�2012-1.5-12-000-1003-016� and the RFFI grant "Spectral and
asymptotic methods for studying of the differential operators" No
11-01-00458.

\end{document}